\documentclass[preprint,10pt]{elsarticle}
\usepackage{hyperref}
\pdfoutput=1
\usepackage{enumerate}
\usepackage{graphics}
\usepackage{epsfig}
\usepackage{amsmath}
\usepackage{amsfonts}
\usepackage{amssymb}
\usepackage{latexsym}

\newtheorem{theorem}{Theorem}[section]
\newtheorem{lemma}{Lemma}[section]

\newcommand{\ds}{\displaystyle}

\newcommand{\eps}{\varepsilon}

\usepackage{amssymb}

\begin{document}

\begin{frontmatter}


\title{ A parameter uniform essentially first order convergent numerical method for a parabolic singularly perturbed differential
 equation of reaction-diffusion type with initial and Robin boundary conditions}

\author{R.Ishwariya}
\address{Department of Mathematics, Bishop Heber College, Tiruchirappalli, Tamil Nadu, India.}
\ead{ishrosey@gmail.com}

\author{J.J.H.Miller}
\address{Institute for Numerical Computation and Analysis, Dublin, Ireland.}
\ead{jm@incaireland.org}

\author{S.Valarmathi\corref{mycorrespondingauthor}}
\address{Department of Mathematics, Bishop Heber College, Tiruchirappalli, Tamil Nadu, India.}
\cortext[mycorrespondingauthor]{Corresponding author}
\ead{valarmathi07@gmail.com}

\begin{abstract}
In this paper, a class of linear parabolic singularly perturbed second order differential equations of reaction-diffusion type with
initial and Robin boundary conditions is considered. The solution $ u $ of this equation is smooth, whereas $\dfrac{\partial u}{\partial x}$ exhibits parabolic boundary layers. 
A numerical method composed of a classical finite difference scheme on a piecewise
uniform Shishkin mesh is suggested. This method is proved to be first order convergent in time and essentially first order convergent in the space variable in
the maximum norm uniformly in the perturbation parameters.
\end{abstract}

\end{frontmatter}
\section{Introduction}

\qquad A differential equation in which small parameters multiply the highest order derivative and some or none of the lower order derivatives
is known as a singularly perturbed differential equation. In this paper, a class of linear parabolic singularly perturbed second order differential equation of reaction-diffusion type with
initial and Robin boundary conditions is considered. \\

For a general introduction to parameter-uniform numerical methods for singular perturbation problems, see \cite{11}, \cite{12}, \cite{18} and \cite{19}. In \cite{13}, a Dirichlet boundary value problem for a linear parabolic singularly perturbed differential equation is studied and a numerical method 
comprising of a standard finite difference operator on a fitted piecewise uniform mesh is considered and it is proved to be uniform with respect to the small
parameter in the maximum norm. In \cite{14}, a boundary-value problem for a singularly perturbed parabolic PDE with convection is considered on an interval
in the case of the singularly perturbed Robin boundary condition is considered and using a defect correction technique, an ε-uniformly convergent schemes of
high-order time-accuracy is constructed. The efficiency of the new defect-correction schemes is confirmed by numerical experiments. In \cite{15}, a 
one-dimensional steady-state convection dominated convection-diffusion problem with Robin boundary conditions is considered and the numerical solutions 
obtained using an upwind finite difference scheme on Shishkin meshes are uniformly convergent with respect to the diffusion cofficient.\\

Consider the following parabolic initial-boundary value problem for a singularly perturbed linear second order
differential equation
\begin{equation}
 \displaystyle\frac{\partial u}{\partial t}(x,t)-\eps\frac{\partial^2 u}{\partial x^2}(x,t)+ a(x,t) u(x,t)= f(x,t), \text{ on } \Omega, \label{31g}
 \end{equation}with
\begin{equation}\label{31aag} 
 \begin{array}{l}
  u(0,t)-\dfrac{\partial  u}{\partial x}(0,t)= \phi_L(t),\;\;
 u(1,t)+\dfrac{\partial  u}{\partial x}(1,t)= \phi_R(t),\;\; 0\leq t \leq T,\\
 \hspace{3.0cm}u(x,0)= \phi_B(x),\;\;0\leq x \leq 1,
 \end{array}
  \end{equation}
where $\Omega = \{(x,t): 0<x<1,\; 0<t\leq T\},\;\; \bar{\Omega}=\Omega \cup \Gamma,\;\; \Gamma = \Gamma_L\cup\Gamma_B\cup\Gamma_R$ with
$\Gamma_L=\{(0,t) :0\leq t\leq T \},\;
\Gamma_R=\{(1,t) :0\leq t\leq T \}\text{ and } \Gamma_B=\{(x,0) :0 < x < 1 \}.$ 

\noindent The problem \eqref{31g}, \eqref{31aag} can also be written in the operator form
\begin{equation*} L u =  f\; \text{ on }\; \Omega,\label{31mg}\end{equation*}
\begin{equation*}\beta_0 u(0,t)= \phi_L(t),\;\;\beta_1 u(1,t)=  \phi_R(t),\;\;
 u(x,0)= \phi_B(x),\label{32mg}\end{equation*}
where the operators $ L, \beta_0, \beta_1 $ are 
defined by \begin{equation*} L = \displaystyle \frac{\partial}{\partial t}-\varepsilon\frac{\partial^2}{\partial x^2}+ a,\;
\beta_0=I-\dfrac{\partial}{\partial x},\;\beta_1=I+\dfrac{\partial}{\partial x} \end{equation*}
where $I$ is the identity operator. The reduced problem corresponding to \eqref{31g}, \eqref{31aag} is defined by
\begin{equation}  \displaystyle\frac{\partial u_0}{\partial t}+ a u_0= f,\;  \text{ on } \;\Omega,\;\;
 u_0= u\; \text{ on }\; \Gamma_B. \label{32g}\end{equation}
The problem \eqref{31g}, \eqref{31aag} is said to be singularly perturbed in the following sense.\\
The solution $u$ of \eqref{31g}, \eqref{31aag} is expected to exhibit weak twin layers of width $O(\sqrt{\eps})$ at $x=0$ and $x=1.$ 
\section{ Solution of the continuous problem}\label{c3s12}
\qquad Standard theoretical results on the existence of the solution of \eqref{31g}, \eqref{31aag} are stated, without proof, in this section. 
See \cite{16} and \cite{17} for more details.
For all $(x,t) \in \bar\Omega,$ it is assumed that $a(x,t)$ satisfies the condition 
 \begin{equation}\displaystyle0<\alpha<a(x,t).\label{34g}\end{equation}
\noindent Sufficient conditions for the existence, uniqueness and regularity of solution of \eqref{31g}, \eqref{31aag} are given in the following theorem.
\begin{theorem} \label{3stg} Assume that $ a \text{ and } f$ are sufficiently smooth.  Also assume that $\phi_L \in C^{2}(\Gamma_L),\;\phi_B \in C^{5}(\Gamma_B)$,\; $  \phi_R \in C^{2}(\Gamma_R)$  and that the following compatibility conditions are fulfilled at the corners $(0,0)$ and $(1,0)$ of\; $\Gamma.$
\begin{equation}\label{36g} \phi_B(0) = \phi_L (0)+\dfrac{d  \phi_B}{dx}(0)\;\;\; \text{and}\;\;\; \phi_B(1) =   \phi_R (0)-\dfrac{d  \phi_B}{dx}(1),\end{equation}
\begin{equation}\label{37g}
\begin{array}{lcl}
\ds\frac{d \phi_L}{dt}(0)&=& - \varepsilon  \ds\frac{d^3 \phi_B}{dx^3}(0)+ \varepsilon  \ds\frac{d^2 \phi_B}{dx^2}(0)+ a(0,0) \ds\frac{d \phi_B}{dx}(0) 
-[ a(0,0)-\dfrac{\partial  a}{\partial x}(0,0)] \phi_B(0)\\&&+  f (0,0)-\dfrac{\partial   f}{\partial x}(0,0),\\
\ds\frac{d   \phi_R}{dt}(0) &=&  \varepsilon  \ds\frac{d^3 \phi_B}{dx^3}(1)+ \varepsilon  \ds\frac{d^2 \phi_B}{dx^2}(1)- a(1,0) \ds\frac{d \phi_B}{dx}(1) 
-[ a(1,0)+\dfrac{\partial  a}{\partial x}(1,0)] \phi_B(1)\\&&+  f (1,0)+\dfrac{\partial   f}{\partial x}(1,0),
\end{array}
\end{equation}
and
\begin{equation}\label{38g}
\begin{array}{lcl}
\dfrac{d^2  \phi_L}{dt^2}(0)&=&- \varepsilon ^2\dfrac{d^5 \phi_B}{d x^5}(0)+ \varepsilon ^2\dfrac{d^4  \phi_B}{dx^4}(0) + 2 \varepsilon  a(0,0)\dfrac{d^3  \phi_B}{dx^3}(0)
+[-2 \varepsilon   a(0,0)\\&&+4 \varepsilon \dfrac{\partial  a}{\partial x}(0,0)]\dfrac{d^2  \phi_B}{dx^2}(0)
+[-2 \varepsilon \dfrac{\partial  a}{\partial x}(0,0)+3 \varepsilon \dfrac{\partial^2  a}{\partial x^2}(0,0)- a^2(0,0)\\&&+\dfrac{\partial  a}{\partial t}(0,0)]
\dfrac{d  \phi_B}{dx}(0)+[- \varepsilon \dfrac{\partial^2  a}{\partial x^2}(0,0)+ a^2(0,0)
-\dfrac{\partial  a}{\partial t}(0,0)+ \varepsilon \dfrac{\partial^3  a}{\partial x^3}(0,0)\\&&-2 a(0,0)\dfrac{\partial  a}{\partial x}(0,0)
+\dfrac{\partial^2  a}{\partial x \partial t}(0,0)]  \phi_B(0)+[- a(0,0)+\dfrac{\partial  a}{\partial x}(0,0)] 
 f(0,0)\\&&+\dfrac{\partial   f}{\partial t}(0,0)- \varepsilon \dfrac{\partial^3   f}{\partial x^3}(0,0)+ \varepsilon \dfrac{\partial^2   f}{\partial x^2}(0,0)+ a(0,0)\dfrac{\partial   f}{\partial x}(0,0)
-\dfrac{\partial^2   f}{\partial x \partial t}(0,0),
\end{array}
\end{equation}
\begin{equation}\label{39g}
\begin{array}{lcl}
\dfrac{d^2  \phi_R}{dt^2}(0)&=& \varepsilon ^2\dfrac{d^5 \phi_B}{d x^5}(1)+ \varepsilon ^2\dfrac{d^4  \phi_B}{dx^4}(1) - 2 \varepsilon  a(1,0)\dfrac{d^3\phi_B}{dx^3}(1)
+[-2 \varepsilon   a(1,0)\\&&-4 \varepsilon \dfrac{\partial  a}{\partial x}(1,0)]\dfrac{d^2  \phi_B}{dx^2}(1)
+[-2 \varepsilon \dfrac{\partial  a}{\partial x}(1,0)-3 \varepsilon \dfrac{\partial^2  a}{\partial x^2}(1,0)+ a^2(1,0)\\&&-\dfrac{\partial  a}{\partial t}(1,0)]\dfrac{d  \phi_B}{dx}(1)
+[- \varepsilon \dfrac{\partial^2  a}{\partial x^2}(1,0)+ a^2(1,0)
-\dfrac{\partial  a}{\partial t}(1,0)- \varepsilon \dfrac{\partial^3  a}{\partial x^3}(1,0)\\&&+2 a(1,0)\dfrac{\partial  a}{\partial x}(1,0)
-\dfrac{\partial^2  a}{\partial x \partial t}(1,0)]  \phi_B(1)+[- a(1,0)-\dfrac{\partial  a}{\partial x}(1,0)]
 f(1,0)\\&&+\dfrac{\partial   f}{\partial t}(1,0)+ \varepsilon \dfrac{\partial^3   f}{\partial x^3}(1,0)+ \varepsilon \dfrac{\partial^2   f}{\partial x^2}(1,0)
- a(1,0)\dfrac{\partial   f}{\partial x}(1,0)+\dfrac{\partial^2   f}{\partial x \partial t}(1,0).
\end{array}
\end{equation}
\end{theorem}
Then there exists a unique solution $ u$ of \eqref{31g}, \eqref{31aag} satisfying $ u \in C_\lambda ^{(4)} (\bar\Omega)$.

\section{ Analytical results}\label{c3s13}
\qquad The operator $ L$ satisfies the following maximum principle:
\begin{lemma}\label{3maxg}
Let the assumptions \eqref{34g} - \eqref{39g} hold. Let $  \psi$ be any function in the domain of $ L$ such that 
$  \beta_0\psi (0,t)\geq 0,\;   \beta_1\psi (1,t)\geq 0,\; \psi (x,0)\geq 0.$ Then\; $ L\psi (x,t) \geq  0$\; on\; $\Omega$ \;implies that
$\psi (x,t) \geq  0$\; on\; $\bar{\Omega}.$
\end{lemma}
\textbf{Proof}. Let $ x^*, t^*$ be such that $\psi (x^*,t^*) = \displaystyle \min_{\bar\Omega} \psi (x,t)$ and assume that
the lemma is false. Then $\psi (x^*, t^*) < 0 .$ For $x^*=0,$
$  \beta_0\psi(0,t^*)=\psi(0,t^*)-\dfrac{\partial \psi }{\partial x}(0,t^*)<0,\;\;$for $x^*=1, 
  \beta_1\psi(1,t^*)=\psi(1,t^*)+\dfrac{\partial\psi}{\partial x}(1,t^*)<0\;\;$and for $ t^*=0,\;\; \psi (x^*,0)<0,$ contradicting the hypotheses. 
Therefore, $(x^* , t^*) \notin \Gamma.$ Let $(x^*,t^*)\in \Omega.$ Then $\dfrac{\partial \psi }{\partial t}(x^*,t^*)= 0\text{ and } \dfrac{\partial^2 \psi }
{\partial x^2}(x^*,t^*)\geq 0$ lead to
\begin{equation*} L  \psi(x^*,t^*)= \dfrac{\partial \psi }{\partial t}(x^*,t^*)-\varepsilon \dfrac{\partial^2 \psi }
{\partial x^2}(x^*,t^*)+a(x^*,t^*) \psi (x^*,t^*)<0,\end{equation*}
which contradicts the assumption and proves the result for $ L.$
\begin{lemma}\label{3srg} Let the assumptions \eqref{34g} - \eqref{39g} hold. If $\psi $ is any function in the domain of $ L,$ then, \;for each\;
$(x,t)\in\bar\Omega,$\begin{equation*} |\psi(x,t)| \le\;
\max\displaystyle\left\{\parallel   \beta_0  \psi(0,t)\parallel,\parallel   \beta_1  \psi(1,t)\parallel,\parallel   \psi(x,0)\parallel, 
\dfrac{1}{\alpha}\parallel  L \psi\parallel\right\}.\end{equation*}
\end{lemma}
\textbf{Proof}. Define the two functions\begin{equation*}\begin{array}{rcl} \theta^{\pm} (x, t) = \max\displaystyle\left\{\parallel   \beta_0  \psi(0,t)\parallel,\parallel   \beta_1  \psi(1,t)\parallel,\parallel   \psi(x,0)\parallel, 
\dfrac{1}{\alpha}\parallel  L\psi\parallel\right\}  \pm \;\;\psi (x,t), \;\; (x,t)\in \bar\Omega.\end{array}\end{equation*}It is not hard to verify that $  \beta_0{\theta}^{\pm}(0,t)\geq  0,\;  \beta_1{\theta}^{\pm}(1,t)\geq  0,\;{\theta}^{\pm}(x,0)\geq  0$ and
$ L{\theta}^{\pm}\geq 0$ on $\Omega.$ 
It follows from Lemma \ref{3maxg} that ${\theta}^{\pm} \geq  0$ on $\bar{\Omega}$ as required.\\
A standard estimate of the solution $u$ of the problem \eqref{31g}, \eqref{31aag} and its derivatives is contained in the following lemma.
\begin{lemma}\label{3udg}Let the assumptions \eqref{34g} - \eqref{39g} hold and let $ u$ be the solution of the problem \eqref{31g}, \eqref{31aag}.
Then,\; for all $(x, t) \in \bar{\Omega}$
\begin{equation*}
\begin{array}{lcl}
|u(x,t)| &\leq& C( \parallel \phi_L(t) \parallel+\parallel \phi_R(t) \parallel+\parallel \phi_B(x) \parallel + \parallel f \parallel ),\\
|\dfrac{\partial^l  u}{\partial t^l}(x,t)| &\leq& C( \parallel u \parallel +\ds\sum_{q=0}^{l} \parallel \dfrac{\partial^q f}{\partial t^q} \parallel ),\;\; l=1,2,\\
|\dfrac{\partial  u}{\partial x}(x,t)| &\leq& C ( \parallel  u \parallel + \parallel \dfrac{\partial f}{\partial x} \parallel ),\\
|\dfrac{\partial^l  u}{\partial x^l}(x,t)| &\leq& C\eps^{\frac{-(l-1)}{2}}( \parallel  u \parallel + \parallel  f \parallel + \parallel \dfrac{\partial f}{\partial t} \parallel + \parallel \dfrac{\partial f}{\partial x} \parallel + \parallel \dfrac{\partial^2 f}{\partial x\partial t} \parallel ), \;\; l=2,3,\\
|\dfrac{\partial^4  u}{\partial x^4}(x,t)| &\leq& C \eps^{\frac{-3}{2}} ( \parallel  u \parallel + \parallel f \parallel +  \parallel \dfrac{\partial f}{\partial t} \parallel + \parallel \dfrac{\partial^2   f}{\partial t^2} \parallel + \parallel \dfrac{\partial f}{\partial x} \parallel + \parallel \dfrac{\partial^2 f}{\partial x\partial t} \parallel + \parallel \dfrac{\partial^3 f}{\partial x\partial t^2} \parallel \\&&+
\eps^{\frac{1}{2}} \parallel \dfrac{\partial^{2} f}{\partial x^{2}} \parallel ), \\
|\dfrac{\partial^2  u}{\partial x \partial t}(x,t)| &\leq& C( \parallel  u \parallel + \parallel f \parallel + \parallel \dfrac{\partial f}{\partial t} \parallel + \parallel \dfrac{\partial   f}{\partial x} \parallel + \parallel \dfrac{\partial^2   f}{\partial x\partial t} \parallel ),\\
|\dfrac{\partial^3  u}{\partial x^2 \partial t}(x,t)| &\leq& C\eps^{\frac{-1}{2}}( \parallel  u \parallel + \parallel f \parallel + \parallel \dfrac{\partial f}{\partial t} \parallel + \parallel \dfrac{\partial^2 f}{\partial t^2} \parallel + \parallel \dfrac{\partial   f}{\partial x} \parallel + \parallel \dfrac{\partial^2   f}{\partial x\partial t} \parallel + \parallel \dfrac{\partial^3   f}{\partial x\partial t^2} \parallel ).
\end{array}
\end{equation*}
\end{lemma}
\textbf{Proof}. The bound on $ u$ is an immediate consequence of Lemma \ref{3srg}.\\ Differentiating \eqref{31g} partially with respect to $`t`$ once and
twice respectively, and applying Lemma \ref{3maxg}, the bounds on $\dfrac{\partial  u}{\partial t}$ and $\dfrac{\partial^2  {u}}{\partial t^2}$
respectively are derived. Now, differentiating \eqref{31g} partially with respect to $`x` $ once, gives
\begin{equation} \dfrac{\partial^2 u}{\partial x\partial t} - \varepsilon \dfrac{\partial^3 u}{\partial x^3}+a\dfrac{\partial  u}
{\partial x}=\dfrac{\partial   f}{\partial x}-\dfrac{\partial a}{\partial x} u,\label{310g} \end{equation}
and from the initial and boundary conditions, we derive
\begin{equation}
\begin{array}{lcl}
\dfrac{\partial  u}{\partial x}(0,t)= u(0,t)- \phi_L(t),\;\;
\dfrac{\partial  u}{\partial x}(1,t)= \phi_R(t)- u(1,t),\;\;
\dfrac{\partial  u}{\partial x}(x,0)=\dfrac{\partial  \phi_B(x)}{\partial x}=\eta(x).\label{311g}
 \end{array}\end{equation}
Denoting $\dfrac{\partial  u}{\partial x}$ by $ z $ in \eqref{310g} and \eqref{311g}, we get,
\begin{equation}
\dfrac{\partial z }{\partial t}-\varepsilon \frac{\partial^2 z }{\partial x^2}+ a z = h  ,  \label{312g}\end{equation}
\begin{equation}                                                                                     
 z (0,t)= u(0,t)-  \phi_L(t),\;\;  z (1,t)=  \phi_R(t)- u(1,t),\;\;  z (x,0)=\eta (x),\label{313g}\end{equation}
where $ h  =\dfrac{\partial  f}{\partial x}-\dfrac{\partial a}{\partial x} u.$\\
This problem \eqref{312g}, \eqref{313g} is similar to the problem in \cite{13}. Now, using the stability result in \cite{13}, the bound on $ z $
or $\dfrac{\partial u}{\partial x}$ is determined. Thus,\begin{center} $|\dfrac{\partial  u}{\partial x}(x,t)|\leq C (\parallel  u\parallel+\parallel \dfrac{\partial   f}{\partial x}\parallel).$\end{center}
On differentiating (\ref{312g}) partially with respect to $`t`$ once and twice respectively, and applying the stability in \cite{13}, the following bounds on 
 $\dfrac{\partial  z }{\partial t}$ or $\dfrac{\partial^2  {u}}{\partial x \partial t}$ and $\dfrac{\partial^2  {z}}{\partial t^2}$ respectively are derived 
 \begin{equation*}
 \begin{array}{lcl}
 |\dfrac{\partial z}{\partial t}(x,t)| \leq C (\parallel  z  \parallel+\parallel h\parallel+\parallel\dfrac{\partial h}{\partial t} \parallel),\;\;
|\dfrac{\partial^2 z}{\partial t^2}(x,t)|&\leq C (\parallel  z  \parallel+\parallel h\parallel+\parallel\dfrac{\partial  h}{\partial t} \parallel
+\parallel\dfrac{\partial^2  h}{\partial t^2} \parallel).
\end{array}
\end{equation*}
To bound $\dfrac{\partial z}{\partial x },$ for each $(x,t)$, consider an interval $I = [a, a +\sqrt{\varepsilon}],\; a \geq 0$
such that $x \in I.$ Then for some $y$ such that $a < y < a + \sqrt{\varepsilon}$\; and \;$t \in (0,T],$
\begin{center}
 $\dfrac{\partial z}{\partial x}(y,t)=\dfrac{z(a+\sqrt{\varepsilon},t)-z(a,t)}{\sqrt{\varepsilon}}.$\end{center}
  Therefore, 
  \begin{equation}
 |\dfrac{\partial z}{\partial x}(y,t)| \leq C \varepsilon^{\frac{-1}{2}}\parallel  z \parallel.\label{314g}
\end{equation}
Then, for any $x\in I,$\begin{center} $\dfrac{\partial z}{\partial x}(x,t) = \dfrac{\partial z}{\partial x}(y,t)+
\displaystyle\int^x_y\dfrac{\partial^2 z(s,t)}{\partial x^2}ds$\end{center}
\begin{center} $\dfrac{\partial z}{\partial x}(x,t) = \dfrac{\partial z}{\partial x}(y,t)+\varepsilon^{-1}
\displaystyle\int^x_y\left(\dfrac{\partial z(s,t)}{\partial t}-h(s,t)+a(s,t)z(s,t)\right)ds.$\end{center}
Therefore,\begin{center}$|\dfrac{\partial z}{\partial x}(x,t)| \leq |\dfrac{\partial z}{\partial x}(y,t)|+C\varepsilon^{-1}
\displaystyle\int^x_y(\parallel  z \parallel+\parallel h\parallel+\parallel\dfrac{\partial  h}{\partial t} \parallel)ds.$\end{center}
Using (\ref{314g}) in the above equation, yields
\begin{center}$|\dfrac{\partial z}{\partial x}(x,t)| \leq C\varepsilon^{\frac{-1}{2}}(\parallel  z  \parallel+\parallel h\parallel+\parallel\dfrac{\partial  h}{\partial t} \parallel).$\end{center}
\begin{center}\text{i.e. } 
$|\dfrac{\partial^2  u}{\partial x^2}(x,t)| \leq C\varepsilon^{\frac{-1}{2}}(\parallel  u \parallel+\parallel  f\parallel+\parallel\dfrac{\partial   f}{\partial t} \parallel
+\parallel\dfrac{\partial   f}{\partial x} \parallel+\parallel\dfrac{\partial^2   f}{\partial x\partial t} \parallel).$\end{center}
Rearranging the terms in (\ref{312g}), we get
\begin{center}
 $|\dfrac{\partial^2 z}{\partial x^2}(x,t)| \leq C\varepsilon^{-1}(\parallel  z  \parallel+\parallel h\parallel+\parallel\dfrac{\partial h}{\partial t} \parallel).$
\end{center}
\begin{center}\text{i.e. }
$|\dfrac{\partial^3  u}{\partial x^3}(x,t)| \leq C\varepsilon^{-1}(\parallel  u \parallel+\parallel  f\parallel+\parallel\dfrac{\partial   f}{\partial t} \parallel
+\parallel\dfrac{\partial f}{\partial x} \parallel+\parallel\dfrac{\partial^2   f}{\partial x\partial t} \parallel).$\end{center}
Following the steps  similar to those used to bound  $\dfrac{\partial z}{\partial x},$ the bound of the mixed derivative $\dfrac
{\partial^2 z}{\partial x \partial t}$ or $\dfrac{\partial^3 u}{\partial x^2 \partial t}$ is also derived.\\
Differentiating (\ref{312g}) once partially with respect to $`x`$ and rearranging the equation, the bound on $\dfrac{\partial^3  z }{\partial x^3}$ or 
$\dfrac{\partial^4  u}{\partial x^4}$ follows.\\\\
\noindent The Shishkin decomposition of the solution $ u$ of the problem \eqref{31g}, \eqref{31aag} is
\begin{equation*}\label{315ag} u= v+ w, \end{equation*}
where the smooth component $v$ of the solution $u$ satisfies
\begin{equation}\label{315bg}
 L{v}={f}, \text{ on } \Omega 
 \end{equation}
with
 \begin{equation}\label{315cg} 
\beta_0v(0,t)=\beta_0u_0(0,t),\;\;\beta_1v(1,t)= \beta_1u_0(1,t),\;\; v(x,0)=u_0(x,0)\end{equation}
and the singular component $w$ of the solution $u$ satisfies
 \begin{equation}
  L{w}={0}, \text{ on } \Omega \label{316g}
 \end{equation}
with
 \begin{equation}\label{316hg} 
\begin{array}{rcl}
\beta_0{w}(0,t)=\beta_0({u}-{v})(0,t),\;\;
\beta_1{w}(1,t)=\beta_1({u}-{v})(1,t),\;\;
 w(x,0)=0.
 \end{array}
 \end{equation}
 \noindent Bounds on the smooth component $v$ of $u$ and its derivatives are contained in
 \begin{lemma}\label{3vdg} Let the assumptions \eqref{34g} - \eqref{39g} hold. Then there exists a constant $C,$ such that, for each $(x,t) \in \bar \Omega,$ 
  \begin{equation*}
\begin{array}{rclc}
|\dfrac{\partial^l v}{\partial t^l}(x,t)|&\leq& C, \;\; l=0,1,2,\;\;
 &|\dfrac{\partial^l v}{\partial x^l}(x,t)|\leq C, \;\; l=1,2,3\\
 |\dfrac{\partial^4 v}{\partial x^4}(x,t)|&\leq& C \eps^{-1/2},\;\;
&|\dfrac{\partial^{l+1} v}{\partial x^l\partial t}(x,t)|\leq C,\;\; l=1,2.
\end{array}
\end{equation*}
 \end{lemma}
\noindent \textbf{Proof}. 
The smooth component $ v$ is subjected to further decomposition
\begin{equation}\label{315aa}
 v=u_{0}+ \eps v_1.
 \end{equation}
The component $v_{1}$ satisfies the following equation:
 \begin{equation}\label{317gg}
\dfrac{\partial v_1}{\partial t}- \eps\dfrac{\partial^2 v_1}{\partial x^2}+av_1=\dfrac{\partial^2u_{0}}{\partial x^2}
 \end{equation} 
with
 \begin{equation}\label{317hg}
 (v_1-\dfrac{\partial v_1}{\partial x})(0,t)=0,\;\;(v_1+\dfrac{\partial v_1}{\partial x})(1,t)=0,\;\; v_1(x,0)=0,
 \end{equation}
 where $u_0$ is the solution of the reduced problem \eqref{32g}.\\
 From the expressions \eqref{317gg}, \eqref{317hg} and using Lemma \eqref{3udg}, it is found that for $k=1,2,3,4,\;\;l=0,1,2,\;\;m=1,2$
\begin{equation}\label{317ig}
\begin{array}{c}
 |\dfrac{\partial^l v_{1}}{\partial t^l}(x,t)|\leq C,\;|\dfrac{\partial^k v_{1}}{\partial x^k}(x,t)|\leq C(1+\eps^{-(k-1)/2}),\;
 |\dfrac{\partial^{m+1} v_{1}}{\partial x^m\partial t}(x,t)|\leq C(1+\eps^{-(m-1)/2}).
\end{array}
\end{equation}
From \eqref{315aa} and \eqref{317ig}, the following bounds hold:
\begin{equation*}
\begin{array}{rclc}
|\dfrac{\partial^l v}{\partial t^l}(x,t)|&\leq& C, \;\; l=0,1,2,\;\;
 &|\dfrac{\partial^l v}{\partial x^l}(x,t)|\leq C, \;\; l=1,2,3\\
 |\dfrac{\partial^4 v}{\partial x^4}(x,t)|&\leq& C \eps^{-1/2},\;\;
&|\dfrac{\partial^{l+1} v}{\partial x^l\partial t}(x,t)|\leq C,\;\; l=1,2.
\end{array}
\end{equation*}\\
\noindent The layer functions $B^{L},\; B^{R},\; B,$ associated with the solution $ u,$ are
defined on $\bar\Omega$ by
$B^{L}(x) = e^{-x\sqrt{\alpha/\varepsilon}},\;B^{R}(x) =
B^{L}(1-x),\;B(x) = B^{L}(x)+B^{R}(x).$\\
The following elementary properties of these layer functions, for all $0\leq x \le y \leq 1,$ should be noted:\\
$B(x)=B(1-x),\;
B^{L}(x)> B^{L}(y), \;0\;<\;B^{L}(x)\leq\;1,\;
B^{R}(x)< B^{R}(y), \;0\;<\;B^{R}(x)\leq\;1,\;
B(x)\; \text{\;is monotonically decreasing for increasing}\;\; x \in [0,\frac{1}{2}],\\
B(x) \;\text{\;is monotonically increasing for increasing}\;\; x \in [\frac{1}{2},1],\\
B(x) \leq 2B^L(x)\; \text{for}\; x \in
[0,\frac{1}{2}],\;
 \;B(x) \leq 2B^R(x) \;\text{for}\; x
\in [\frac{1}{2},1],\;
B^L(2\frac{\sqrt\varepsilon}{\sqrt\alpha}\ln N)=N^{-2}.$\\\\
\noindent Bounds on the singular component $ w $ of $ u$ and its derivatives are contained in
\begin{lemma}\label{3wdg} Let the assumptions \eqref{34g} - \eqref{39g} hold. Then there exists a constant $C,$ such that, for each $(x,t) \in \bar \Omega,$ 
\begin{equation*}\label{336g} \begin{array} {lcll}
|\dfrac{\partial^l w }{\partial t^l}(x,t)| \le C B(x),\;\text{ for} \;\; l=0,1,2, \;\;
&&|\dfrac{\partial^l w }{\partial x^l}(x,t)| \le C \dfrac{B(x)}{\eps^{\frac{l-1}{2}}},\;\text{ for}\;\; l=1,2,\\
|\dfrac{\partial^3 w }{\partial x^3}(x,t)| \le C \dfrac{B(x)}{\eps},\;\;
&&|\dfrac{\partial^4 w }{\partial x^4}(x,t)| \le C \dfrac{B(x)}{\eps^{\frac{3}{2}}}.
\end{array} \end{equation*}
\end{lemma}
\noindent \textbf{Proof}. 
To derive the bound of $w$, define two functions
\begin{center}${\psi}^\pm(x,t)=Ce^{\alpha t}B (x)\;\pm\;w (x,t)$. \end{center}
For a proper choice of $C,$\; $  \beta_0  \psi^\pm(0,t)\geq 0,\;  \beta_1  \psi^\pm(1,t)\geq 0$\;and\;${ \psi^\pm}(x,0)\geq 0$ and for\;$(x,t)\in \Omega$, 
\[\begin{array}{lcl}
 L \psi^\pm(x,t) &=& C \alpha e^{\alpha t} B (x)-C\eps\frac{\alpha}{\eps}e^{\alpha t} B (x) + Cae^{\alpha t}B (x)\;\ge 0.
\end{array}
\]
By Lemma \ref{3maxg}, $\psi ^{\pm} \geq  0$ on $\bar\Omega$ and it follows that
\begin{equation*}\label{31414g}
 |w (x,t)|\leq Ce^{\alpha t}B (x)\;\;\;\text{or}\;\;\; |w (x,t)|\leq CB (x).
\end{equation*}
Differentiating \eqref{316g} partially with respect to
$`t`$ once and twice, and using Lemma \ref{3maxg}, it is not hard to see that
 \begin{equation*}
|\frac{\partial w }{\partial t}(x,t)|\leq CB (x),\;\;|\frac{\partial^2 w }{\partial t^2}(x,t)|\leq CB (x).
\end{equation*}
\noindent Differentiating \eqref{316g} with respect to $`x` $ once,
\begin{equation} \dfrac{\partial^2 w }{\partial x\partial t} - \eps\dfrac{\partial^3 w }{\partial x^3}+ 
a\dfrac{\partial  w }{\partial x}=-\dfrac{ \partial a}{\partial x} w.\label{338g} \end{equation}
And from the initial and boundary conditions,
\begin{equation}\label{339g}
\begin{array}{lcl}
 \dfrac{\partial  w }{\partial x}(0,t)=  w (0,t) - \beta_0  (u-v) (0,t),\;\;
\dfrac{\partial  w }{\partial x}(1,t)=   \beta_1  (u-v) (1,t) - w (1,t)  ,\;\;
\dfrac{\partial  w }{\partial x}(x,0)= 0.\end{array}
 \end{equation}
Denoting $\dfrac{\partial w}{\partial x}$ by $ \tilde{z}$ in \eqref{338g} and \eqref{339g}, yields
\begin{equation}
\dfrac{\partial{\tilde{z}}}{\partial t}-\eps\frac{\partial^2{\tilde{z}}}{\partial x^2}+ a{\tilde{z}}
=-\dfrac{ \partial a}{\partial x} w ,  \label{340g}\vspace{-0.5cm} \end{equation}
\begin{equation}                                                                                     
\tilde{z}(0,t)= w (0,t) - \beta_0(u-v) (0,t),\;\; \tilde{z}(1,t)= \beta_1(u-v)(1,t) - w (1,t),\;\; \tilde{z}(x,0)=0.\label{341g}\end{equation}
This problem \eqref{340g}, \eqref{341g} is similar to the problem in \cite{13}.  Now, using the stability result in \cite{13}, 
the bound on $\tilde{z}$ or $\dfrac{\partial w}{\partial x}$ is determined. Thus,\begin{equation*}\label{342g} |\dfrac{\partial w }{\partial x}(x,t)|\leq C B (x)\end{equation*}
Differentiating (\ref{340g}) partially with respect to $`t`$ once and twice respectively, and applying the stability in \cite{13}, the following bounds on 
 $\dfrac{\partial  \tilde{z} }{\partial t}$ or $\dfrac{\partial^2 w}{\partial x\partial t}$ and $\dfrac{\partial^2  \tilde{z}}{\partial t^2}$ respectively are derived. 
 \begin{equation*}
 \begin{array}{c}
 |\dfrac{\partial \tilde z}{\partial t}(x,t)| \leq C B(x),\;\;|\dfrac{\partial^2 \tilde z}{\partial t^2}|\leq C B(x).
\end{array}
\end{equation*}
To bound $\dfrac{\partial^2 w }{\partial x^2}$ or $\dfrac{\partial \tilde{z}}{\partial x },$ consider an interval $I = [c, c +\sqrt{\varepsilon}],\; c \geq 0$
such that $x \in I.$ Then for some $y$ such that $c < y < c + \sqrt{\varepsilon}$\; and \;$t \in (0,T],$
\begin{equation*}
 \dfrac{\partial \tilde{z}}{\partial x}(y,t)=\dfrac{\tilde{z}(c+\sqrt{\varepsilon},t)-\tilde{z}(c,t)}{\sqrt{\varepsilon}}.\end{equation*} 
 Therefore,
 \begin{equation*}
 |\dfrac{\partial \tilde{z}}{\partial x}(y,t)| \leq \dfrac{C}{\sqrt{\varepsilon}} (B (c+\sqrt{\eps})+B (c)).
\end{equation*}
Hence, \begin{equation} 
\begin{array}{lcl}
|\dfrac{\partial \tilde{z}}{\partial x}(y,t)|& \leq& \dfrac{C}{\sqrt{\varepsilon}} (B^L(c+\sqrt{\eps})+B^R(c+\sqrt{\eps})+B^L(c)+B^R(c))\\
&\leq &\dfrac{C}{\sqrt{\varepsilon}} (B^L(c)+B^R(c+\sqrt{\eps})).\label{342ag}
\end{array}\end{equation}
Then, for any $x\in I,$ such that $y < \eta <x $
\begin{equation*}\label{343g}\begin{array}{lcl} \dfrac{\partial \tilde{z}}{\partial x}(x,t) &=& \dfrac{\partial \tilde{z}}{\partial x}(y,t)+
\displaystyle\int_y^x\dfrac{\partial^2 \tilde{z}}{\partial x^2}(\eta,t)d\eta\\
&=&\dfrac{\partial \tilde{z}}{\partial x}(y,t)+\eps^{-1}\displaystyle\int_y^x(\dfrac{\partial{\tilde{z}}}{\partial t}+ a{\tilde{z}}+\dfrac{ \partial a}
{\partial x}w)(\eta,t)d\eta\end{array}\end{equation*}
By using the bounds for $\tilde{z},\; \dfrac{\partial \tilde{z}}{\partial t},\; w$ and \eqref{342ag} in the above equation yields
\[\begin{array}{lcl}
|\dfrac{\partial \tilde{z}}{\partial x}(x,t)|  &\leq& C \eps^{\frac{-1}{2}}(B^L(c)+B^R(c+\sqrt{\eps}))+C \eps^{\frac{-1}{2}}(B^L(\eta)+B^R(\eta))\\
 &=& C \eps^{\frac{-1}{2}}(B ^L(x)\dfrac{B ^L(c)}{B ^L(x)} + B ^R(x)\dfrac{B ^R(c+\sqrt{\eps})}{B ^R(x)})\\
 &\leq& C \eps^{\frac{-1}{2}}B (x).
 \end{array}\]
Therefore,
\begin{equation*} 
\begin{array}{lcl}
 |\dfrac{\partial \tilde{z}}{\partial x}(x,t)|\leq \dfrac{C}{ \sqrt{\eps}}B (x)\text{ or }
 |\dfrac{\partial^2 w }{\partial x^2}(x,t)|\leq \dfrac{C}{ \sqrt{\eps}}B (x). \label{1422g} 
\end{array}
\end{equation*}
Rearranging the equation \eqref{340g} yields
\begin{equation*}
\eps\frac{\partial^2{\tilde{z}}}{\partial x^2}
= \dfrac{\partial{\tilde{z}}}{\partial t}+ a\tilde{z}+\dfrac{ \partial a}{\partial x} w.\end{equation*}
Using the bounds of $w,$ $\tilde z$ and $\dfrac{\partial \tilde z}{\partial t}$ in the above equation, the following bound holds.
\begin{equation*}\label{342bg}\begin{split}|\dfrac{\partial^2\tilde{z}}{\partial x^2}(x,t)|\leq C \eps^{-1}B(x), \;\;\text{or\;\;}
|\dfrac{\partial^3w}{\partial x^3}(x,t)|\leq C \eps^{-1}B(x).
\end{split}\end{equation*}
Differentiating \eqref{340g} with respect to $`t`$ once and 
following a similar procedure to bound $\dfrac{\partial \tilde z}{\partial x},$ the bound of the mixed derivative $\dfrac{\partial^2 \tilde z}{\partial x\partial t}$
 or $\dfrac{\partial^3 w}{\partial x^2\partial t}$ is derived.
\begin{equation*} 
\begin{array}{lcl}
 |\dfrac{\partial^3 w }{\partial x^2\partial t}(x,t)|\leq C\varepsilon^{\frac{-1}{2}} B(x).
\end{array}
\end{equation*}
On differentiating \eqref{340g} with respect to $x$ and rearranging yields
\[|\eps \frac{\partial^3 \tilde{z}}{\partial x^3}(x,t)| \leq C \eps^{-1/2}B(x),\;\;
\text{or\;\;}|\eps \frac{\partial^4 w }{\partial x^4}(x,t)| \leq C \eps^{-1/2}B(x).\]

\section{The Shishkin mesh}\label{c3s14}
A piecewise uniform Shishkin mesh is now
constructed. Let $\Omega^M_t=\{t_k \}_{k=1}^{M},\;\Omega^N_x=\{x_j
\}_{j=1}^{N-1},\;\bar{\Omega}^M_t=\{t_k
\}_{k=0}^{M},\;\bar{\Omega}^N_x=\{x_j
\}_{j=0}^{N},\;\Omega^{N,M}=\Omega^N_x \times \Omega^M_t,\;
 \bar{\Omega}^{N,M}=\bar{\Omega}^N_x \times \bar{\Omega}^M_t \;\text{ and }\;\Gamma^{N,M}=\Gamma \cap \bar{\Omega}^{N,M}.$\\
 \noindent The mesh $\bar{\Omega}^M_t$ is chosen to be a uniform mesh with $M$
sub-intervals on $[0,T]$. The mesh $\bar{\Omega}^N_x$ is a
piecewise-uniform Shishkin mesh with $N$ mesh intervals. The interval $[0,1]$ is subdivided into 
$3$ sub-intervals 
\[[0,\sigma]\cup(\sigma,1-\sigma]\cup(1-\sigma,1]\]
where
\begin{equation*}\label{3151g}\sigma = \min\displaystyle\left\{\frac{1}{4},\; 2\frac{\sqrt\varepsilon}{\sqrt\alpha}\ln
N\right\}\end{equation*}
Thus, on the sub-interval $(\sigma,1-\sigma],$ a uniform mesh with
$\;\frac{N}{2}$ mesh-points is placed and on the subintervals $(0,\sigma]$ and $(1-\sigma,1],\;$a uniform mesh with
$\;\frac{N}{4}$ mesh-points is placed. \\
In particular, when the parameter $ \sigma $ takes on its left-hand value, the Shishkin mesh $\bar{\Omega}_x^N$ becomes a classical uniform mesh.
In practice, it is convenient to take $N=4k,\; k \geq 2.$\\
 \indent From the above construction of $\bar{\Omega}_x^N$ it is clear that the transition points $\{\sigma,1-\sigma\}$
are the only points at which the mesh size can change and that it does not necessarily
change at each of these points. The following notations are introduced: $h_j= x_j-x_{j-1} ;\; h_{j+1} = x_{j+1}- x_j,\; J=\{x_j=\sigma,1-\sigma : h_{j+1}\neq h_{j}\}.$
For each point $x_j$ in the sub-intervals $(0,\sigma]$ and $(1-\sigma,1],\; x_j-x_{j-1}=4N^{-1}\sigma $ and for $x_j$ in the sub-interval $(\sigma,1-\sigma],$ 
$x_j-x_{j-1}=2N^{-1}(1-2\sigma).$ 

 The construction of $\bar\Omega_x^N$ as a piecewise uniform Shishkin mesh on $[0,1]$ leads to a piecewise uniform Shishkin mesh $\bar\Omega^{N,M}$ on
 $[0,1]\times[0,1]$ 
by considering the cartesian product of the discrete set $\bar\Omega_x^N=\{x_j\}_{j=0}^{N}$ with $\bar\Omega_t^M=\{t_j\}_{j=0}^{M}.$\\
Thus $\bar\Omega^{N,M}$ is a piecewise uniform Shishkin grid with $NM$ mesh elements.\\

\section{The discrete problem}\label{c3s15}
In this section, a classical finite difference operator with an
appropriate Shishkin mesh is used to construct a numerical method
for the problem \eqref{31g}, \eqref{31aag}, which is shown later to be first order parameter-uniform convergent in time and essentially first order
parameter-uniform convergent in the space variable.

The discrete initial-boundary value problem is now defined by the finite difference scheme on the Shishkin mesh $\bar\Omega^{N,M},$ defined in the previous section. 
\begin{equation}\label{3161g} 
 D^-_t{U}(x_j,t_k)-\eps\delta^2_x{U}(x_j,t_k) +a(x_j,t_k){U}(x_j,t_k)= f(x_j,t_k)\;\; \text{ on }\;\; \Omega^{N,M},
 \end{equation}
with
 \begin{equation}\label{3161ag}
 \begin{array}{lcl}
  U(0,t_k) - D^+_x  U(0,t_k) =  \phi_L(t_k),\;\;
  U(1,t_k) + D^-_x  U(1,t_k) =  \phi_R(t_k),\\
\hspace{4.0cm} U(x_j,0) =  \phi_B(x_j).
\end{array}
 \end{equation}
The problem \eqref{3161g}, \eqref{3161ag} can also be written in the operator form
\[ L^{N,M} {U}= f \;\; \text{ on }
\;\; \Omega^{N,M},\]
\begin{equation*}\label{3161bg}
 \begin{array}{lcl}
   \beta_0^{N,M}{U}(0,t_k) =  \phi_L(t_k),\;\;
   \beta_1^{N,M}{U}(1,t_k) =  \phi_R(t_k),\;\;
 U(x_j,0) =  \phi_B(x_j),\end{array}
 \end{equation*}
where $ L^{N,M} =  D^-_t-\eps\delta^2_x+a,\;\;\beta_0^{N,M}=I-D_x^+,\;\;\beta_1^{N,M}=I+D_x^-$.\\
The following discrete results are analogous to those for the
continuous case.
\begin{lemma}\label{3dmaxg} Let the assumptions \eqref{34g} - \eqref{39g} hold.
Then, for any mesh function $ \Psi  $, the inequalities 
$  \beta_0^{N,M} \Psi  (0,t_k) \ge  0  ,\;  \beta_1^{N,M} \Psi  (1,t_k) \ge  0  ,\; \Psi  (x_j,0) \ge  0  $  and $ L^{N,M}
 \Psi  \;\ge\; 0  $ on\; $\Omega^{N,M}$ imply that $
{\Psi}\ge  0  $ on $\bar{\Omega}^{N,M}.$
\end{lemma}
\noindent \textbf{Proof}. Let $ j^*, k^*$ be such that
$ \Psi   (x_{j^{*}}, t_{k^{*}})=\ds\min_{j,k} \Psi   (x_j,t_k)$
and assume that the lemma is false. Then $ \Psi   (x_{j^{*}},
t_{k^{*}})<0$ . From the hypotheses,\;
$ \Psi    (x_{j^*},t_{k^{*}}) -  \Psi    (x_{j^*},t_{k^{*}-1}) < 0$, $ \Psi    (x_{j^*},t_{k^{*}})- \Psi   (x_{j^*-1},t_{k^{*}})<
0, \;  \Psi   
(x_{j^*+1},t_{k^{*}})- \Psi   (x_{j^*},t_{k^{*}})> 0.$ Hence, 
$\;D_t^- \Psi   (x_{j^*}, t_{k^*}) < 0$\;and\;$\delta^2_x \Psi   (x_{j^*},t_{k^{*}})\;\geq\;0.$ It follows that, for $(x_{j^*},t_{k^*})\in\Omega^{N,M},$ 
\[\begin{array}{lcl} L^{N,M} \Psi (x_{j^*},t_{k^{*}})\;&=&
\;D_t^- \Psi   (x_{j^*}, t_{k^*})-\eps\delta^2_x \Psi   (x_{j^*},t_{k^{*}})+a(x_{j^*},t_{k^{*}})\Psi(x_{j^*},t_{k^{*}})\;<0,\end{array}\]
which is a contradiction.\\
If $x_{j^*}=0,$ then 
$ \beta_0^{N,M}\Psi(0,t_{k^{*}})= \Psi(0,t_{k^{*}})\;-\;D^+_x \Psi   (0,t_{k^{*}})<0,$ a contradiction. Therefore, $x_{j^*}\neq0,$ and for a similar 
reason $x_{j^*}\neq1.$ For $t_{k^*}=0,\;\Psi(x_{j^*},0)<0,$ which is a contradiction. Therefore, $t_{k^*}\neq0.$  
Hence the result.\\
\noindent An immediate consequence of this is the following discrete stability
result.
\begin{lemma}\label{3dsrg} Let the assumptions \eqref{34g} - \eqref{39g} hold.
Then, for any mesh function $ \Psi   $ defined on $\bar\Omega^{N,M},$
\[|{\Psi}(x_j,t_k)|\;\le\;\max\left\{ \parallel   \beta_0^{N,M} \Psi  (0,t_k) \parallel ,\; \parallel \beta_1^{N,M} \Psi(1,t_k) \parallel ,\; \parallel  \Psi(x_j,0) \parallel ,\; \frac{1}{\alpha} \parallel 
 L^{N,M} \Psi   \parallel \right\}. \]
\end{lemma}
\noindent \textbf{Proof}. Define the two mesh functions
\begin{equation*}\begin{array}{rcl}{\Theta}^{\pm}(x_j,t_k)=\max\{ \parallel   \beta_0^{N,M} \Psi  (0,t_k) \parallel ,\; \parallel   \beta_1^{N,M} \Psi  (1,t_k) \parallel ,\; \parallel \; \Psi \;  (x_j,0) \parallel ,\; \frac{1}{\alpha} \parallel 
 L^{N,M} \Psi   \parallel \}\\\pm  \Psi  (x_j,t_k),\;(x_j,t_k)\in \bar\Omega^{N,M}.\end{array}\end{equation*}
 It is not hard to verify that
 $  \beta_0^{N,M}{\Theta}^{\pm}(0,t_k)\geq  0  ,\;  \beta_1^{N,M}{\Theta}^{\pm}(1,t_k)\geq  0  ,\;{\Theta}^{\pm}(x_j,0)\geq  0  $ and
$ L^{N,M}{\Theta}^{\pm}\geq  0  $ on $\Omega^{N,M}$. It follows from Lemma \ref{3dmaxg} that ${\Theta}^{\pm}\geq  0  $ on
$\bar{\Omega}^{N,M}$.\\
\noindent The following comparison principle will be used in the proof of the error estimate.
\begin{lemma}\label{3compg} Assume that the mesh functions ${\Phi}$ and  ${Z}$ satisfy\;\;
$|\beta_0^{N,M} Z(0,t_k)| \leq\\ \beta_0^{N,M}\Phi(0,t_k),\;| \beta_1^{N,M} Z(1,t_k)| \leq  \beta_1^{N,M}\Phi(1,t_k),\;|Z(x_j,0)| \leq \Phi(x_j,0)\text{ and }
| L^{N,M}{Z}| \leq  L^{N,M} {\Phi}\\ \text{ on }\Omega^{N,M}.$ Then 
$|Z| \leq \Phi\text{ on } \bar\Omega^{N,M}.$
\end{lemma}
\noindent \textbf{Proof}. Define the two mesh functions $ \Psi  ^{\pm}$ by \;\;$ \Psi  ^{\pm}={\Phi} \pm {Z}.$\\
Then, $ \Psi^{\pm}$
satisfy
 $ \beta_0^{N,M}\Psi^{\pm}(0,t_k) \geq 0,\; \beta_1^{N,M}\Psi^{\pm}(1,t_k) \geq 0,\; \Psi^{\pm}(x_j,0) \geq 0,\text{ and }
 L^{N,M} \Psi  ^{\pm} \geq 0\text{ on }\Omega^{N,M}.$ The required result follows from the Lemma \ref{3dmaxg}.
\section{The local truncation error}\label{c3s16}
From Lemma \ref{3dsrg}, it is seen that in order to bound the error
$ U - u$, it suffices to bound 
$ \beta_0^{N,M}(U-u)(0,t_k),\;\beta_1^{N,M}(U-u)(1,t_k),\;(U-u)(x_j,0),\text{ and } L^{N,M}( U - u)$. Note that, for $(x_j, t_k) \in \Omega ^{N,M},$
$ L^{N,M}( U - u)= L^{N,M}U - L^{N,M}u= f- L^{N,M} u
= L u- L^{N,M} u
=( L- L^{N,M}) u.$\\
It follows that $ L^{N,M}( U - u)=(\dfrac{\partial}{\partial
t}-D^-_t) u-\eps(\dfrac{\partial^2}{\partial
x^2}-\delta^2_x) u.$\\
Let $ V , W $ be the discrete analogues of
$ v , w $ respectively, given by
\begin{equation*}
L^{N,M} V = f \text{ on } \Omega^{N,M},
\end{equation*}
\begin{equation*}
\beta_0^{N,M}V(0,t_k)=\beta_0v(0,t_k),\; \beta_1^{N,M}V(1,t_k)=\beta_1v(1,t_k),\; V(x_j,0)=v(x_j,0), 
\end{equation*}
\begin{equation*}
L^{N,M}W = 0 \text{ on } \Omega^{N,M},
\end{equation*}
\begin{equation*}
\beta_0^{N,M}W(0,t_k)=\beta_0w(0,t_k),\; \beta_1^{N,M}W(1,t_k)=\beta_1w(1,t_k),\; W(x_j,0)=w(x_j,0),
\end{equation*}
where $v$ and $ w$ are the solutions of  \eqref{315bg}, \eqref{315cg} and \eqref{316g}, \eqref{316hg} respectively. Further, 
\begin{equation*}|\beta_0^{N,M}(V-v)(0,t_k)|=|(\frac{\partial}{\partial x}-D^+_x)v(0,t_k)|,\;
|\beta_1^{N,M}(V-v)(1,t_k)|=|(D^-_x-\frac{\partial}{\partial x})v(1,t_k)|,\end{equation*}
\begin{equation*}|\beta_0^{N,M}(W-w)(0,t_k)|=|(\frac{\partial}{\partial x}-D^+_x)w(0,t_k)|,\;
|\beta_1^{N,M}(W-w)(1,t_k)|=|(D^-_x-\frac{\partial}{\partial x})w(1,t_k)|,\end{equation*}
\begin{equation*}| L^{N,M}( V   - v )(x_j,t_k)|=|((\frac{\partial}{\partial
t}-D^-_t) v -\eps(\frac{\partial^2}{\partial x^2}-\delta^2_x) v)(x_j,t_k)| ,\end{equation*}
\begin{equation*}|L^{N,M}( W - w )(x_j,t_k)|=|((\frac{\partial}{\partial
t}-D^-_t) w -\eps(\frac{\partial^2}{\partial
x^2}-\delta^2_x) w)(x_j,t_k)|.\end{equation*}
\noindent The local truncation error of the smooth and singular components can be treated separately. Note that, for any smooth function $ \psi $ and for each
$(x_j,t_k)\in \Omega^{N,M}$, the following
distinct estimates of the local truncation error hold:
\begin{equation}\label{3174g}
|(\frac{\partial}{\partial t}-D^-_t)\psi(x_j,t_k)|\;\le\;
C(t_k-t_{k-1})\max_{s\;\in\;[t_{k-1},\;t_k]}|\frac{\partial^2\psi}{\partial
t^2}(x_j,s)|,
\end{equation}
\begin{equation}\label{3175ag}
|(\frac{\partial}{\partial x}-D^-_x)\psi(x_j,t_k)|\;\le\;
C(x_{j}-x_{j-1})\max_{s\;\in\;[x_{j-1},\;x_{j}]}|\frac{\partial^2\psi}{\partial
x^2}(s, t_k)|,
\end{equation}
\begin{equation}\label{3175g}
|(\frac{\partial}{\partial x}-D^+_x)\psi(x_j,t_k)|\;\le\;
C(x_{j+1}-x_{j})\max_{s\;\in\;[x_{j},\;x_{j+1}]}|\frac{\partial^2\psi}{\partial
x^2}(s, t_k)|,
\end{equation}
\begin{equation}\label{3176g}
\hspace{-2.5cm}|(\frac{\partial^2}{\partial x^2}-\delta^2_x)\psi(x_j,t_k)|\;\le\;
C\max_{s\;\in\;I_j}|\frac{\partial^2\psi}{\partial x^2}(s,t_k)|,
\end{equation}
\begin{equation}\label{3177g}
\hspace{-.5cm}|(\frac{\partial^2}{\partial
x^2}-\delta^2_x)\psi(x_j,t_k)|\;\le\;C(x_{j+1}-x_{j-1}) \max_{s\;\in\;I_j}|\frac{\partial^3\psi}{\partial x^3}(s,t_k)|.
\end{equation}
Here $I_j=[x_{j-1}, x_{j+1}]$.
\section{Error estimate}\label{c3s17}
The proof of the theorem on the error estimate is broken into two parts. First, a theorem concerning the error in the smooth component is established.
Then the error in the singular component is estimated.\\
Define the barrier function through
\begin{equation*}\label{3182g}{\Phi}(x_j,t_k)=C[M^{-1}+2N^{-1}\ln N+
\sigma N^{-1}\ln N \theta(x_j,t_k)],\end{equation*} where $C$ is sufficiently large and 
$\theta$ is a piecewise linear polynomial on $\bar\Omega$ for each $x_j=\sigma\in J$ defined by
\begin{equation*}\theta(x,t)=
\left\{ \begin{array}{l}\;\; \dfrac{x}{ \sigma    }, \;\; 0 \leq x \leq  \sigma,\\
\;\; 1, \;\;  \sigma     < x < 1- \sigma,\\
\;\; \dfrac{1-x}{ \sigma    }, \;\; 1- \sigma     \leq x \leq 1. \end{array}\right .\end{equation*}\\
Also note that,
\begin{equation}\label{3181g} L^{N,M} \theta(x_j,t_k) \ge
\left\{ \begin{array}{l}\;\;
\alpha  \theta   (x_j, t_k), \;\;  \text{ if} \;\;x_j \notin J \\
\;\; \alpha+\dfrac{\eps N}{\sigma}, \;\;  \text{ if} \; x_j \in J. \end{array}\right .\end{equation}
Then, on $\Omega^{N,M}$,  ${\Phi}$ satisfies\qquad
$0 \leq \Phi (x_j,t_k) \leq
C(M^{-1}+N^{-1}\ln N).$\\
Also, 
\begin{equation*}\label{3184ag}
\begin{split}
\beta_0^{N,M}\Phi(0,t_k) \ge C(M^{-1}+N^{-1}\ln N),\;\;
\beta_1^{N,M}\Phi(1,t_k) \ge C(M^{-1}+N^{-1}\ln N).
 \end{split} 
\end{equation*}
\noindent For $x_j\notin J$, using (\ref{3181g}), it is not hard to see that,
\begin{equation*}\label{3184g}
L^{N,M}{\Phi}(x_j,t_k) \ge C(M^{-1}+N^{-1}\ln N)\end{equation*}
and, for $x_j \in J$, using (\ref{3181g}), it is not hard to see that, 
\begin{equation*}\label{3185g}
L^{N,M}{\Phi}(x_j,t_k) \ge C(M^{-1}+N^{-1}\ln N).\end{equation*}
The following theorem gives the estimate of the error in the smooth component $V$.
 \begin{theorem}\label{3smootherrorthmg} Let the assumptions \eqref{34g} - \eqref{39g} hold. Let $v$ denote the smooth component of the 
solution of the problem \eqref{31g}, \eqref{31aag} and $V$ denote the smooth component of the
solution of the problem \eqref{3161g}, \eqref{3161ag}.  Then
\begin{equation*}\;\;  \parallel  V   - v  \parallel  \leq C(M^{-1}+N^{-1}\ln N).\end{equation*}\end{theorem}
\noindent \textbf{Proof}. 
From the expression \eqref{3175g},
\begin{equation}\label{3179abg} 
\begin{array}{rcl} 
|\beta_0^{N,M}(V-v)(0,t_k)| &\leq&C(x_{1}-x_{0})\ds\max_{s\;\in\;[x_{0},\;x_{1}]}|\dfrac{\partial^2 v}{\partial
x^2}(s, t_k)|\\
&\leq&C N^{-1}.
 \end{array}
 \end{equation}
From the expression \eqref{3175ag},
\begin{equation}\label{3179acg} 
\begin{array}{rcl} 
|\beta_1^{N,M}(V- v)(1,t_k)| &\leq&C(x_{N}-x_{N-1})\ds\max_{s\;\in\;[x_{N-1},\;x_{N}]}|\dfrac{\partial^2 v}{\partial
x^2}(s, t_k)|\\
&\leq&C N^{-1}.
 \end{array}
 \end{equation}
 Also note that, the expressions \eqref{3184ag}, \eqref{3179abg} and \eqref{3179acg} yield
\begin{equation}\label{3567g} 
\begin{array}{lcl}
|\beta_0^{N,M}( V- v)(0,t_k)| \leq \beta_0^{N,M}\Phi(0,t_k),\;
|\beta_1^{N,M}( V- v)(1,t_k)| \leq \beta_1^{N,M}\Phi(1,t_k),\\
\hspace{4.0cm}|( V- v)(x_j,0)| \leq \Phi(x_j,0),
 \end{array}
 \end{equation}
For each mesh point $x_j,$ there are two possibilities: either $x_j \notin J$ or $x_j \in J$.\\
For $x_j \notin J$, from the expressions \eqref{3174g}, \eqref{3177g} and Lemma \ref{3vdg}, it follows that
\begin{equation}\label{3179g}
\begin{array}{lcl}
| L^{N,M}({V}-{v}) (x_j,t_k)| &\leq& |((\dfrac{\partial}{\partial t}-D^-_t)-\eps(\dfrac{\partial^2}{\partial x^2}-\delta^2_x)) v (x_j,t_k)|\\
&\leq& C[M^{-1}+N^{-1}]\\
&\leq& L^{N,M}\Phi(x_j,t_k)
\end{array}\end{equation}
For $x_j \in J$, then $x_j=\sigma \text{ or }1-\sigma.$ Here the argument for $x_j=\sigma$ is given and for
$x_j=1-\sigma$ it is analogous.\\
For $x_j =\sigma,$ from the expressions \eqref{3174g}, \eqref{3177g} and Lemma \ref{3vdg}, it follows that
\begin{equation}\begin{array}{lcl}
| L^{N,M}({V}-{v})(x_j,t_k)|& \leq& |((\dfrac{\partial}{\partial t}-D^-_t)-\eps(\dfrac{\partial^2}{\partial
x^2}-\delta^2_x)) v (x_j,t_k)|\\
&\leq& C[M^{-1}+ N^{-1}\ln N]\\
&\leq& L^{N,M}\Phi(x_j,t_k)\label{3189g}
\end{array}\end{equation}\\
From the expressions \eqref{3567g}, \eqref{3179g}, \eqref{3189g} and the comparison principle, the required result follows.\\
\noindent The following theorem gives the estimate of the error in the singular component $\vec W$.
\begin{theorem} \label{3singularerrorthmg} Let the assumptions \eqref{34g} - \eqref{39g} hold. Let $w$ be the singular component of the
solution of the problem \eqref{31g}, \eqref{31aag} and $ W$ be the singular component of the
solution of the problem \eqref{3161g}, \eqref{3161ag}.  Then
\begin{equation*}\;\;  \parallel  W - w  \parallel  \leq C(M^{-1}+N^{-1}\ln N). \end{equation*}
\end{theorem}
\noindent \textbf{Proof}.
From the expression \eqref{3175g},
\begin{equation}\label{3179bb} 
\begin{array}{rcl} 
|\beta_0^{N,M}(W- w)(0,t_k)| &\leq&C(x_{1}-x_{0})\ds\max_{s\;\in\;[x_{0},\;x_{1}]}|\dfrac{\partial^2 w}{\partial
x^2}(s, t_k)|\\
&\leq&C N^{-1}\ln N,
 \end{array}
 \end{equation}
From the expression \eqref{3175ag},
\begin{equation}\label{3179bc} 
\begin{array}{rcl} 
|\beta_1^{N,M}(W- w)(1,t_k)| &\leq&C(x_{N}-x_{N-1})\ds\max_{s\;\in\;[x_{N-1},\;x_{N}]}|\dfrac{\partial^2 w}{\partial
x^2}(s, t_k)|\\
&\leq&C N^{-1}\ln N.
 \end{array}
 \end{equation}
 Also note that,
\begin{equation}\label{3190cg} 
\begin{array}{lcl}
|\beta_0^{N,M}(W-w)(0,t_k)|\leq \beta_0^{N,M}\Phi(0,t_k),\;
|\beta_1^{N,M}(W-w)(1,t_k)|\leq \beta_1^{N,M}\Phi(1,t_k),\\
\hspace{4.0cm}|(W- w )(x_j,0)|\leq \Phi(x_j,0).
 \end{array}
 \end{equation}
\noindent For $x_j \notin J,$ from the expressions \eqref{3174g}, (\ref{3177g}) and Lemma \ref{3wdg},
it follows that
\begin{equation}\label{3190ag}
\begin{array}{lcl}|L^{N,M}( W   - {w})(x_j,t_k)|&=&
|((\dfrac{\partial}{\partial t}-D^-_t)-\eps(\dfrac{\partial^2}{\partial
x^2}-\delta^2_x)) w (x_j,t_k)|\\ 
& \leq &
C(M^{-1}+(x_{j+1}-x_{j-1})\;\ds\max_{s\;\in\;I_j}{B(s)}) \\
& \leq & C(M^{-1}+N^{-1})\\
&\leq& L^{N,M}\Phi(x_j,t_k)
\end{array}\end{equation}
\noindent For $x_j\in J, $ from the expressions \eqref{3174g}, \eqref{3177g} and Lemma \ref{3wdg}, it follows that
\begin{equation}\label{3190bg}
 \begin{array}{lcl}|L^{N,M}( W   - {w})(x_j,t_k)|&=&
|((\dfrac{\partial}{\partial t}-D^-_t)-\eps(\dfrac{\partial^2}{\partial
x^2}-\delta^2_x)) w (x_j,t_k)|\\ 
& \leq &
C(M^{-1}+(x_{j+1}-x_{j-1})\;\ds\max_{s\;\in\;I_j}{B(s)}) \\
& \leq & C(M^{-1}+N^{-1}\ln N)\\
&\leq& L^{N,M}\Phi(x_j,t_k)
\end{array}\end{equation}
From the expressions \eqref{3190cg}, \eqref{3190ag}, \eqref{3190bg} and the comparison principle, the required result follows.\\
The following theorem gives a parameter uniform bound which is first order in time and essentially
first order in space for the convergence of the discrete solution.
\begin{theorem}\label{3maing} Let the assumptions \eqref{34g} - \eqref{39g} hold. Let $u$ denote the solution of the problem \eqref{31g}, \eqref{31aag} 
and $U$ denote the solution of the problem \eqref{3161g}, \eqref{3161ag}.  Then
\begin{equation*}  \parallel  U - u \parallel  \leq C(M^{-1}+N^{-1}\ln N).\end{equation*}
\end{theorem}
 \textbf{Proof.}
An application of the triangular inequality and the results of
Theorem \ref{3smootherrorthmg} and Theorem \ref{3singularerrorthmg} lead to the required result.
\section{Numerical Illustration}\label{c3s18}
\qquad The numerical method proposed above is illustrated through the example presented in this section. 
The method proposed above is applied to solve the
problem and the parameter-uniform order of convergence and the parameter-uniform error constants are computed. 
To get the order of convergence in the variable $t$ seperately, a Shishkin mesh is considered for $x$ and 
the resulting problem is solved for various uniform meshes with respect to $t$. In order to get the order of convergence in the variable $x$ seperately, a uniform mesh
is considered for $t$ and the resulting problem is solved for various piecewise uniform Shishkin meshes with respect to
$x$. The same two-mesh algorithm found in \cite{12} is applied to get parameter-uniform order of 
convergence and the error constants. The numerical results are presented in Table \ref{t31g} and Table \ref{t32g}.

{\noindent\bf Example  }
Consider the problem 
\[\frac{\partial u}{\partial t} -\eps\frac{\partial^2  u}{\partial x^2} + (1+3t) u=e^{3t}\text{ on }
(0,1) \times (0,1],\]\[ (u-\dfrac{\partial u}{\partial x})(0,t)=1+t^5,\;\;
(u+\dfrac{\partial u}{\partial x})(1,t)=1+t^5,\;\;u(x,0)= 1.\]

For various values of $\eps,$ the maximum errors, the $\eps$- uniform order of convergence and the $\eps$-uniform error constant are computed.
Fixing a Shishkin mesh on $[0,1]$ with $64$ points horizontally, the problem is solved by the method suggested above. 
The order of convergence and the error constant for $u$ are calculated for $t$ using two-mesh algorithm and the results are presented in Table \ref{t31g}. 
A uniform mesh on $[0,1]$ with $256$ points vertically is considered and the order of convergence and the error constant for $u$ in the variable $x$ using two-mesh algorithm are calculated and 
the results are presented in Table \ref{t32g}.

The notations $ D^N,\; p^N $ and $C_{p^{*}}^N $ denote the $\varepsilon$-uniform maximum pointwise two-mesh differences, the $\varepsilon$-uniform order of 
convergence and the $\varepsilon$-uniform error constant respectively and are given by 
\noindent$D^N=\displaystyle\max_{\varepsilon}\; D^N_{\varepsilon}$ where $D^N_{\varepsilon}=\;\parallel U^N_{\varepsilon}- U^{2N}_{\varepsilon} \parallel_{\Omega^N}$, 
$p^N=log_2 \displaystyle\frac{D^N}{D^{2N}}$ and $C_{p^{*}}^N =\displaystyle\frac{D^N N^{p^*}}{1-2^{-p^{*}}}.$ Then the parameter-uniform order of convergence and the error constant are given by 
$p^*=\displaystyle\min_N p^N$ and $C^*_{p^*}=\displaystyle\max_N C^N_{p^*}$ respectively.
\begin{table}[!ht]
\footnotesize
\begin{center}
 \caption{\label{t31g}Values of \;\;$D_{\varepsilon}^N, D^N, p^{N}, p^* \text{ and }C_{p^*}^N\text{ for }\; \alpha = 0.9\;\text{and}\;N=64.$} 
\end{center}
 
\begin{center}
\begin{tabular}{|l|l|l|l|l|}
\hline         
\multicolumn{1}{|c|}{$\eps$} &\multicolumn{4}{|c|}{Number of mesh points $N$}\\
\cline{2-5}& \;\;\; 32 & \;\;\; 64 & \;\;\; 128 & \;\;\; 256\\\hline
\;\;$2^{-6}$ & 0.231E-01   &    0.122E-01   &    0.638E-02  &     0.330E-02\\\hline
\;\;$2^{-8}$ & 0.253E-01   &    0.131E-01   &    0.669E-02   &    0.338E-02\\\hline
\;\;$2^{-10}$ & 0.263E-01   &    0.134E-01   &    0.675E-02   &    0.339E-02\\\hline
\;\;$2^{-12}$ & 0.265E-01    &   0.134E-01   &    0.676E-02  &     0.339E-02  \\\hline
\;\;$2^{-14}$ & 0.266E-01    &   0.134E-01  &     0.676E-02   &    0.339E-02 \\\hline\hline
\; $D^N$ & 0.266E-01    &   0.134E-01   &    0.676E-02   &    0.339E-02\\\hline
 \; $p^N$ & 0.983E+00  &     0.991E+00   &    0.996E+00 & \\\hline
 \; $C_{p^*}^N$ & 0.162E+01  &     0.162E+01  &     0.161E+01  &     0.160E+01\\\hline\hline
\multicolumn{5}{|c|}{Computed $t$-order of $\varepsilon-$uniform convergence, $p^*$ = 0.9827155} \\\hline
\multicolumn{5}{|c|}{Computed $\varepsilon-$uniform error constant, $C_{p^*}^*$ = 1.620163}\\\hline
\end{tabular}
\end{center}
\end{table}
 
\begin{table}[!ht]
\footnotesize
\begin{center}
 \caption{\label{t32g}Values of \;\;$D_{\varepsilon}^N, D^N, p^{N}, p^* \text{ and }C_{p^*}^N\text{ for }\; \alpha = 0.9\;\text{and} \;M=256.$ }
\end{center}

\begin{center}
\begin{tabular}{|l|l|l|l|l|}
\hline         
\multicolumn{1}{|c|}{$\eps$} &\multicolumn{4}{|c|}{Number of mesh points $N$}\\
\cline{2-5}& \;\;\; 32 & \;\;\; 64 & \;\;\; 128 & \;\;\; 256 \\\hline
\;\;$2^{-6}$ &    0.964E-02  &     0.399E-02  &     0.139E-02  &     0.496E-03 \\\hline
\;\;$2^{-8}$ &  0.119E-01    &   0.508E-02   &    0.195E-02    &   0.560E-03 \\\hline
\;\;$2^{-10}$ &  0.117E-01   &    0.617E-02   &    0.258E-02   &    0.843E-03\\\hline
\;\;$2^{-12}$ &  0.537E-02   &    0.298E-02    &   0.155E-02   &    0.698E-03 \\\hline
\;\;$2^{-14}$ &   0.272E-02  &     0.150E-02   &    0.771E-03   &    0.290E-03\\\hline\hline
\; $D^N$ & 0.119E-01   &    0.617E-02   &    0.258E-02   &    0.843E-03 \\\hline
 \; $p^N$ & 0.946E+00   &    0.126E+01  &     0.162E+01 &\\\hline
 \; $C_{p^*}^N$ & 0.655E+00   &    0.655E+00   &    0.528E+00   &    0.332E+00\\\hline\hline
\multicolumn{5}{|c|}{Computed $x$-order of $\varepsilon-$uniform convergence, $p^*$ = 0.9456793 } \\\hline
\multicolumn{5}{|c|}{Computed $\varepsilon-$uniform error constant, $C_{p^*}^*$ = 0.6552203 }\\\hline
\end{tabular}
\end{center}
\end{table}                

\begin{figure}
\begin{minipage}{0.45\textwidth}\caption{\label{f31g}}
\begin{center} The numerical approximation \\of $u$ for $\eps=2^{-14}$ and $M=256$\end{center}
\includegraphics[width=\textwidth]{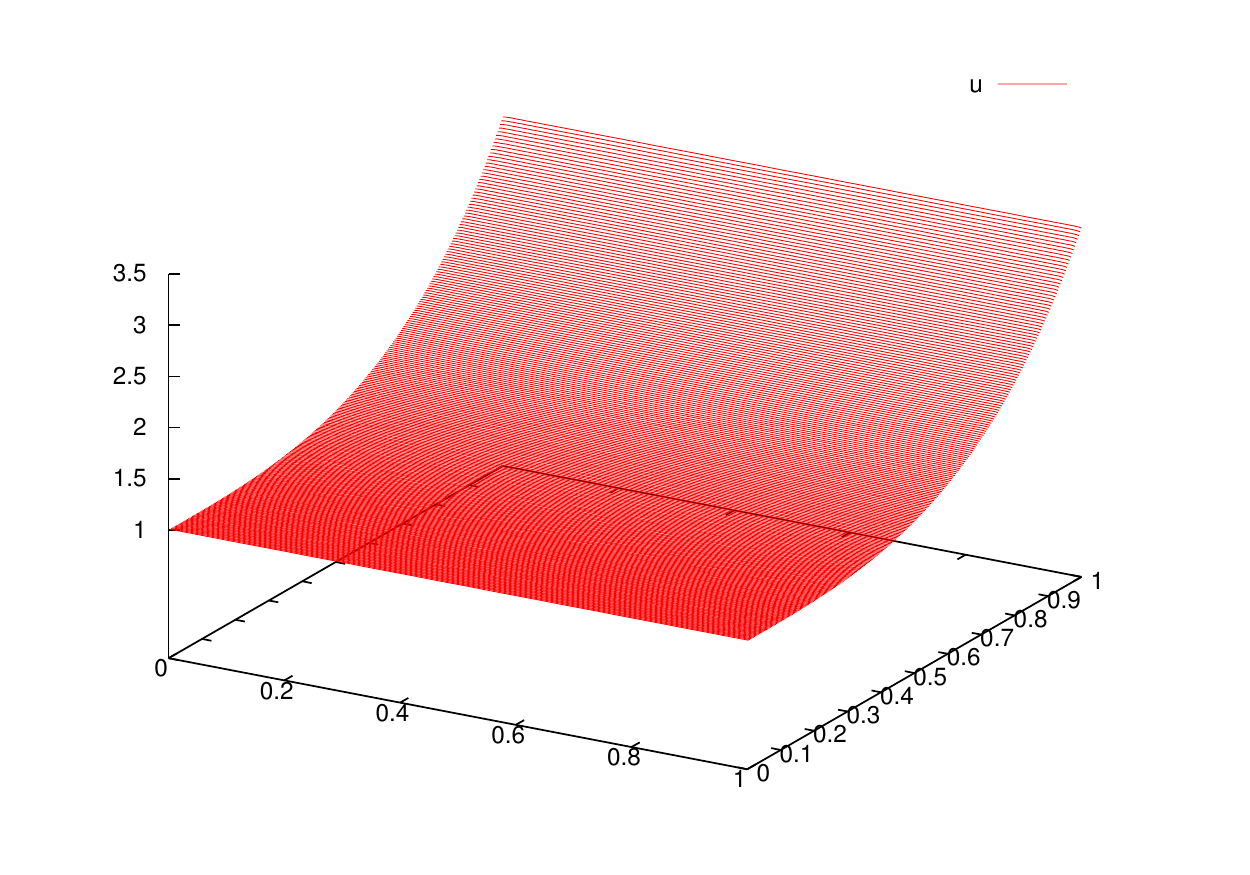} 
\end{minipage} \;\;\;\;
\begin{minipage}{0.45\textwidth}\caption{\label{f32g}}
\begin{center} The numerical approximation \\of $u$ for $\eps=2^{-14}$ and $N=64$\end{center}
\includegraphics[width=\textwidth]{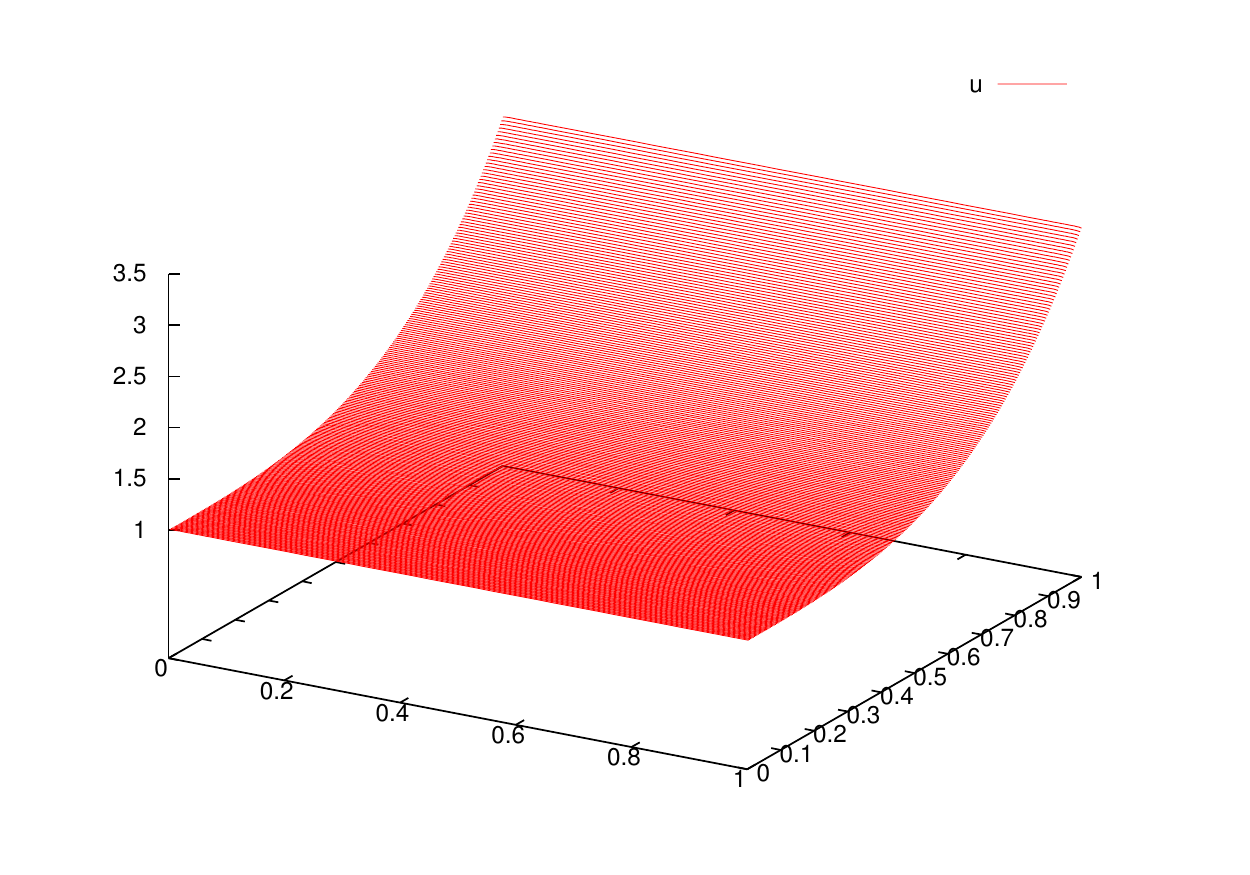} 
\end{minipage} 
\begin{minipage}{0.45\textwidth}\caption{\label{f33g}}
\begin{center}The numerical approximation \\of $\dfrac{\partial u}{\partial x}$ for $\eps=2^{-14}$ and $M=256$\end{center}  
\includegraphics[width=\textwidth]{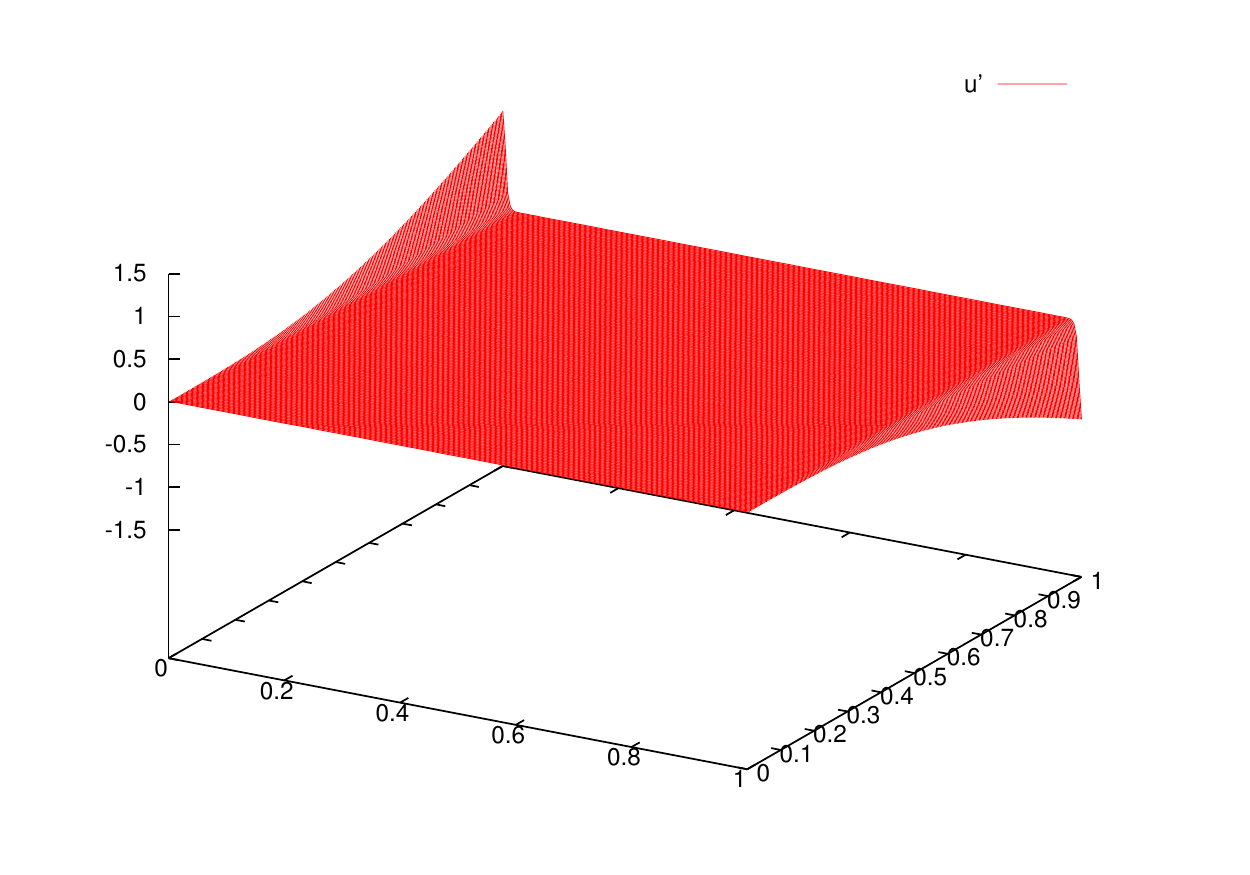} 
\end{minipage}\;\;\;\;
\begin{minipage}{0.45\textwidth}\caption{\label{f34g}}
\begin{center}The numerical approximation \\of $\dfrac{\partial u}{\partial x}$ for $\eps=2^{-14}$ and $N=64$\end{center} 
\includegraphics[width=\textwidth]{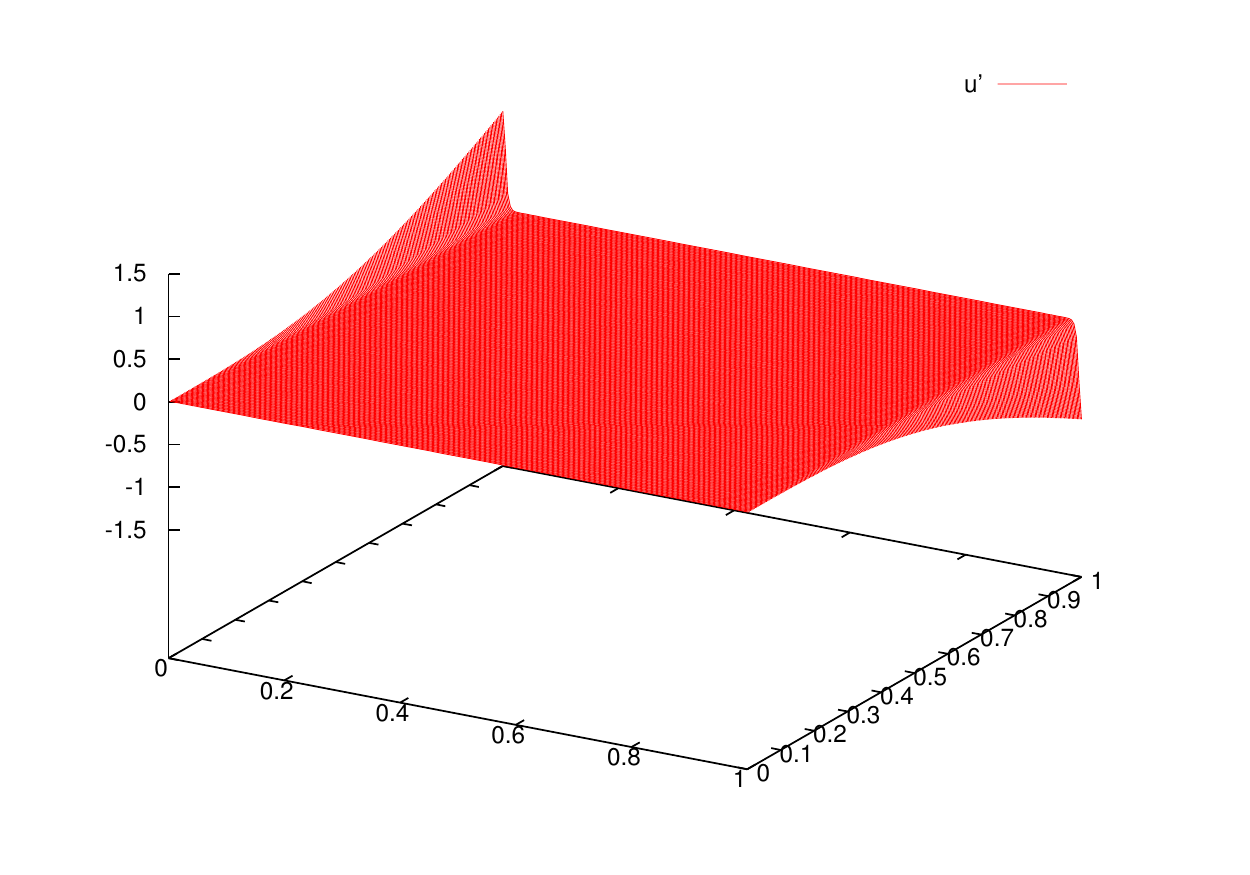} 
\end{minipage}
\end{figure}
It is evident from the Figures \ref{f31g}-\ref{f34g} that the solution $ u$ exhibits no layer whereas the derivative $\dfrac{\partial  u}{\partial x}$ exhibits
parabolic twin boundary layers at $ (0,t) $ and $ (1,t), \; 0 \leq t \leq 1. $ Further, the $t$- order of convergence and the $x$- order of convergence of
the numerical method presented in Table \ref{t31g} and Table \ref{t32g} agree with the theoretical result. 

\section*{Acknowledgment}

The first author sincerely thanks the University Grants Commission, New Delhi, India, for the financial support through the Rajiv Gandhi National
Fellowship to carry out this work.


\end{document}